\documentclass[11pt]{article}   
\usepackage{amssymb,amscd,latexsym}   

\newcommand{\rar}{\rightarrow}
\newcommand{\lar}{\longrightarrow}

\newtheorem{Theorem}{Theorem}[section]

\newtheorem{Corollary}[Theorem]{Corollary}

\newtheorem{Proposition}[Theorem]{Proposition}

\newtheorem{Remark}[Theorem]{Remark}

\newtheorem{Example}[Theorem]{Example}

\newtheorem{definition}[Theorem]{Definition}
\newtheorem{Question}[Theorem]{Question}

\def\demo{\noindent{\em Proof. }}
\newcommand{\Rees}{{\cal R}}
\newcommand{\Spec}{\mbox{\rm Spec}}

\newcommand{\rank}{\mbox{\rm rank }}

\newcommand{\Hom}{\mbox{\rm Hom}}
\newcommand{\Ext}{\mbox{\rm Ext}}

\newcommand{\depth}{\mbox{\rm depth }}

\newcommand{\height}{\mbox{\rm height }}

\newcommand{\Ass}{\mbox{\rm Ass}}

\newcommand{\edim}{\mbox{\rm edim }}
\newcommand{\trdeg}{\mbox{\rm trdeg}}

\def\phi{\varphi}
\def\surjects{\twoheadrightarrow}

\def\dd{{\mathbb D}}
\def\XX{{\bf X}}

\def\TT{{\bf T}}

\def\ass{{\rm Ass}\,}

\def\ecodim{{\rm ecodim}\,}
\def\hht{{\rm height}}
\def\edim{{\rm edim}\,}

\def\spec{{\rm Spec}}

\def\pp{{\mathbb P}}

\def\fm{{\mathfrak m}}
\def\fn{{\mathfrak n}}
\def\cc{{\mathbb C}}

\def\za{\mathbb{S}_{A/k}}
\def\zac{\mathbb{S}_{A/\cc}}

\def\ra{\mathbb{R}_{A/k}}
\def\rac{\mathbb{R}_{A/\cc}}
\def\ba{\mathbb{B}_{A/k}}

\def\QED{\hfill$\Box$}
\def\qed{\QED}

\begin{document}

\title{\LARGE\sc  Tangent Algebras
\vspace{-1mm}}
\footnotetext{AMS 2000 {\it Mathematics Subject
Classification}: Primary 13H10, 14F10; Secondary 14M05, 14M10.\\
\indent {\it Keywords}: K\"ahler differentials, Rees algebra, symmetric
algebra, Zariski tangent algebra.}

\author{
{\normalsize\sc Aron Simis\thanks{Partially supported by CNPq,
Brazil.}}
\vspace{-0.75mm}\\
{\small Departamento de Matem\'atica}\vspace{-1.4mm} \\
{\small Universidade Federal de Pernambuco}\vspace{-1.4mm}\\
{\small 50740-540 Recife, PE, Brazil}\vspace{-1.4mm} \\
{\small e-mail: {aron@dmat.ufpe.br}}\\
\and
{\normalsize\sc Bernd Ulrich\thanks{Partially
 supported by the NSF.}}
\vspace{-0.75mm}\\
{\small Department of Mathematics}\vspace{-1.4mm} \\
{\small Purdue University}\vspace{-1.4mm}\\
{\small West Lafayette, Indiana 47907-1395}\vspace{-1.4mm} \\
{\small e-mail: {ulrich@math.purdue.edu}} \\
 \and
{\normalsize\sc Wolmer V. Vasconcelos\thanks{Partially
 supported by the NSF.}}
\vspace{-0.75mm}\\
{\small Department of Mathematics}\vspace{-1.4mm} \\
{\small Rutgers University}\vspace{-1.4mm}\\
{\small 110 Frelinghuysen Rd}\vspace{-1.4mm}\\
{\small Piscataway, New Jersey 08854-8019}\vspace{-1.4mm} \\
{\small e-mail: {vasconce@math.rutgers.edu}}\vspace{4mm}
}

\maketitle


\begin{abstract}

\noindent We study the Zariski tangent cone $T_X\stackrel{\pi}{\lar}
X$ to an affine variety $X$ and
the closure $\overline{T}_X$ of  $\pi^{-1}({\rm Reg}(X))$ in $T_X$.
We focus on the comparison between $T_X$ and $\overline{T}_X$,
giving sufficient conditions on $X$ in order that
$T_X=\overline{T}_X$. 
One aspect of the results is to understand when this equality takes
place in the presence of the reducedness of the Zariski tangent cone.
Our other interest is to consider conditions on $X$ in order that
$\overline{T}_X$ be normal or/and Cohen--Macaulay, 
and to prove that they are met by several classes of affine varieties including
complete intersection, Cohen--Macaulay codimension two and
Gorenstein codimension three singularities.
In addition, when $X$ is the affine
cone over a smooth
arithmetically normal Calabi--Yau projective variety, we establish 
when $\overline{T}_X$ is also (the affine cone over) an
arithmetically normal Calabi--Yau like (projective) variety. 

\end{abstract}


\newpage
\section{Introduction}

Let $A$ be an affine algebra over a perfect field $k$, and let $\Omega_{A/k}$ be the
$A$-module of K\"ahler $k$-differentials. Classically, the properties of this module
are closely related to the local singularities of $A$, embodying in particular the
well-known Jacobian criterion for the smoothness of $A$. Its sheaf version for an
algebraic variety over a field is fundamental in intersection theory and its
cohomology is a major vehicle for the study of the global geometry of the variety.

In this paper we focus on two basic algebras associated to $\Omega_{A/k}$, where $A$
is essentially of finite type over a perfect field $k$: the Rees algebra $\ra$ of
$\Omega_{A/k}$ and its close predecessor $\za$, the symmetric algebra of
$\Omega_{A/k}$. If $A$ is regular these algebras coincide, but otherwise they may be
quite apart and their respective properties have different impact on the nature of
$A$. In \cite{star} the ring $\za$ was called the
Zariski tangent algebra of $A$ because the closed fibers of the map $\spec (\za)\lar \spec (A)$  are the Zariski
tangent spaces to closed points of $\spec (A)$,  when $A$ is an affine algebra over an algebraically closed field.
Alternatively,  $\spec (\za)$ is the first jet
scheme of $\spec(A)$.
As to $\ra$ it plays the
role of the coordinate ring of a correspondence in biprojective space (see
\cite{elam99}). A variation on $\ra$ is the Rees algebra of the top wedge module of
$\Omega_{A/k}$ which, as was argued in \cite{SSU}, gives a full-fledged algebraic
version of Zak's inequality for the dimension of the Gauss image of a projectively
embedded variety.

If $A$ is reduced, $\Omega_{A/k}$ is generically free and there
is a natural surjection
$$\za\surjects \mathbb{R}_{A/k}.$$
A great deal of the
present work has to do with this map, whose kernel measures
the failure of describing $\ra$ solely in terms of linear equations, or  to use a recent terminology for
modules, the failure of $\Omega_{A/k}$ being of {\it linear
type}.
We also study how the singularities of $A$ are reflected
in the normality, Cohen--Macaulayness and Gorensteiness of $\mathbb{R}_{A/k}$.

Let us now describe the main results in some detail.

\medskip

The paper is divided into two sections. In the first section we focus on the Zariski
tangent algebra $\za$. 
Recall that $\Omega_{A/k}\simeq \dd/\dd^2$, where $\dd$ is the kernel of
the multiplication map $A\otimes_k A\lar A$.
Quite generally, let $S$ be a Noetherian
ring and let ${\bf D}\subset S$ be an ideal such that the symmetric algebra ${\cal S}_A({\bf D}/{\bf D}^2)$ is
torsionfree over $A=S/{\bf D}$. If one assumes,
moreover, that ${\bf D}$ is generically a complete intersection and $A$ is reduced
then ${\cal S}_A({\bf D}/{\bf D}^2)$ is reduced. The converse to this statement,
namely, that ${\cal S}_A({\bf D}/{\bf D}^2)$ is torsionfree if it is reduced, is
known to be false in general. The corresponding question for the
associated graded ring ${\rm gr}_{\bf D}(S)$ was treated in \cite{HuSiVa} and shown
to be affirmative provided $A$ has finite projective dimension over $S$.
Although the diagonal ideal $\dd$ does not have finite projective dimension for singular $A$, our goal
in the first part of the section is to consider this converse in the framework of
$\za={\cal S}_A(\dd/\dd^2)$, trading off the homological restriction on $A$ for a condition on its defining
equations (Theorem~\ref{redtorsionfree} and Corollary~\ref{quadrics}).

In the second part of this section we deal with the same question in the case of algebroid curves. Here, without
any assumptions on the defining equations, we prove that such curves are non-singular provided $\za$ is reduced
(Theorems~\ref{Bergerlike} and \ref{Bergerlike2}).
This can be regarded as an analogue of Berger's conjecture (see \cite{Berger}) in which the reducedness of
$\za$ replaces the torsionfreeness of $\Omega_{A/k}$.

Looking for other natural regularity assumptions on $\Omega_{A/k}$ which imply, in
the spirit of Berger's conjecture, the regularity of $A$, we also prove a version in
which $\Omega_{A/k}$ modulo its torsion is a Cohen--Macaulay $A$-module
(Theorem~\ref{Bergerlite}). Here $A$ is no longer one-dimensional, instead the
standing assumption is that it has embedding codimension $2$ and is sufficiently
well-behaved locally in codimension $2$.

In the last part of the section we study the relationship between
the reflexivity and the normality of $\za$. Quite generally, the
former is known to imply the latter. We will prove that the converse
holds under suitable conditions on the defining equations of $A$ of
order $2$ (Proposition~\ref{normality_quadrics} and
Theorem~\ref{normality_reflexiveness}).

\medskip

The second section is mainly devoted to the behavior of the Rees algebra $\ra$.

Our first fundamental result deals with the case where $A$ is locally everywhere a complete
intersection over a field of characteristic zero and $\ra$ is Cohen--Macaulay
(Theorem 3.1). Notice that, in this case, the module of differentials $\Omega_{A/k}$ has
projective
dimension one.
To better situate the reader, we recall the parallel case of the
symmetric algebra. Quite generally, if $A$ is a local Cohen--Macaulay ring and $E$
is a finite $A$-module of projective dimension, then $E$ satisfies condition $(F_0)$
if and only if the symmetric algebra ${\cal S}_A(E)$ is a complete intersection and $E$
satisfies condition
$(F_1)$ if and only if ${\cal S}_A(E)$ is $A$-torsionfree (see \cite[Proposition 4]{Avramov},
\cite[1.1]{Huneke}, \cite[3.4]{sv}; cf. Section~\ref{RAA} for the notion of condition $(F_t)$).
As the torsionfreeness of  ${\cal S}_A(E)$ gives an isomorphism ${\cal S}_A(E)\simeq
{\cal R}_A(E)$, it
follows that condition $(F_1)$ implies the Rees
algebra ${\cal R}_A(E)$ to be Cohen--Macaulay. However, the converse does not hold in
general, even if
one assumes the preliminary condition $(F_0)$. Nevertheless we are able to prove this converse
in the case where $E$ is the
module of differentials of a complete intersection over a field of characteristic
zero (Theorem 3.1). Our result shows that Cohen--Macaulayness is a rather restrictive property
for
$\ra$. It implies for instance that if $X\subset {\mathbb P}^{2d+1}_k$ is a
$d$-dimensional smooth non-degenerate complete intersection over a field $k$ of
characteristic zero with homogeneous coordinate ring $A$, then $\ra$ is never a
Cohen--Macaulay ring. It would be interesting to have a  geometric explanation of this
phenomenon.

The second kind of structural results concerns the normality of
$\ra$, in which the condition $(F_2)$ will play a predominant role.
Throughout $A$ is a normal domain. To enlarge the picture we
consider yet another algebra associated to an $A$-module $E$,
namely, ${\cal B}_A(E):=\bigoplus_{i\geq 0}(E^i)^{**}$ which could
be dubbed the {\it reflexive closure\/} of the Rees algebra ${\cal
R}_A(E)$ -- note it is a Krull domain with the same divisor class
group as $A$, but may fail to be Noetherian.
This algebra has been studied earlier in \cite{HSV2} (see also
\cite[Chapter 7]{alt}). Here we focus more closely on the case where
$E=\Omega_{A/k}$. First is a basic criterion for the equality ${\cal
R}_A(E)={\cal B}_A(E)$ assuming that ${\cal R}_A(E)$ is normal or
just satisfies Serre's condition $(S_2)$
(Proposition~\ref{normalityandreflexiveness}). The criterion is
given in terms of a new invariant that controls the growth of the
local analytic spreads of $E$. We then move on to an encore of
modules of projective dimension one over Cohen--Macaulay normal
domains.  For such a module $E$ it is known that the condition
$(F_2)$ holds if and only if ${\cal S}_A(E)\simeq {\cal B}_A(E)$
(see \cite[Proposition 4]{Avramov}). In particular, if $E$ satisfies
$(F_2)$ then ${\cal R}_A(E)$ is normal. Surprisingly, the converse
holds as well provided the non-free locus of $E$ be contained in the
singular locus of $A$ (Theorem~\ref{projdim1andnormality}). An
application to $\Omega_{A/k}$ when $A$ is a normal complete
intersection now ensues (Corollary~\ref{CIandnormality}). In fact it
turns out that the inequalities ${\rm edim}\, A_{\mathfrak p}\leq
2\dim A_{\mathfrak p}-2$ locally on the singular locus of $A$ become
equivalent to the powers of $\Omega_{A/k}$ being integrally
closed in the range $1\leq i\leq\, {\rm ecodim}\, A$.

A substantial part of the section is a  push 
 to extend the previous results to
situations other than complete intersections. We have succeeded in the cases where
$A$ is of low dimension or embedding codimension or is sufficiently structured.
Thus, we find satisfactory answers to
 to the
codimension $2$ perfect and codimension $3$ Gorenstein cases
(Propositions~\ref{normcodim2} and \ref{normcodim3}). Naturally the last part
requires several developments
 and the results become  somewhat too technical to
be described in this introduction. We refer the reader to the appropriate places in
the last part of this work.

The section ends with some consequences of the theory so far regarding its relation to 
Calabi--Yau varieties.
We argue that if $A$ is the homogeneous coordinate ring of a Calabi--Yau variety
often $\za$ or $\ra$ are the homogeneous coordinate rings of ``Calabi--Yau like'' varieties.


\section{The symmetric algebra of the module of
differentials}

We look at the symmetric algebra $\za:={\cal S}_A(\Omega_{A/k})$ of the module of
K\"ahler differentials $\Omega_{A/k}$ as an ancestor of the corresponding Rees algebra,
which will be the topic of the next section.

\subsection{Reduced Zariski tangent algebras}

The main result of this part has a curious geometric consequence. Let $A$ be a local
domain with an isolated singularity that is essentially of finite type
over a perfect field.
If the Zariski tangent algebra of $A$ is reduced, but not a domain, then the defining
ideal of $A$ has to contain ``many'' quadrics.
In the language of jet schemes, the same conclusion holds if the first jet scheme of $\spec (A)$
is reduced, but not irreducible.

\begin{Theorem}\label{redtorsionfree}
Let $A$ be a reduced ring essentially of finite type over a perfect field $k$.
Assume that for every non-minimal  ${\mathfrak p}\in \spec (A)$,
$A_{\mathfrak p}\simeq R/I$ where  $(R,{\mathfrak n})$ is a regular local
ring essentially of finite type over $k$ and $I\subset {\mathfrak
n}^2$ is an $R$-ideal satisfying
\[\mu(I+{\mathfrak n}^3/{\mathfrak n}^3)\leq \dim R-1.\]

Then $\za$ is reduced {\rm (}if and{\rm )} only if it is $A$-torsionfree.
\end{Theorem}
\demo The ``if'' statement was explained above.

To show the converse,  let ${\cal T}\subset \za$ denote the
$A$-torsion submodule of $\za$ and suppose that ${\cal T}\neq 0$.
Then there exists
an associated prime of $\za$ contracting to a non-minimal prime of
$A$. We may further take such a non-minimal prime to be
minimal among all non-minimal primes of $A$ that are
contracted from some minimal prime of $\za$. By localizing at this
non-minimal prime, we do not change either the hypotheses or the
conclusion of the statement, so we can reduce the argument to the situation in
which $(A,{\mathfrak m},K)$ is local and ${\mathfrak m}$ is the
contraction of minimal prime of $\za$. Moreover, every minimal prime of $\za$ not
containing ${\mathfrak m}$ must contract to a minimal prime of $A$,
hence contains the torsion ${\cal T}$. From this follows, since
$\za$ is reduced, the crucial relation
\begin{equation}\label{zerointersection}
{\cal T}\cap {\mathfrak m}\za=0.
\end{equation}

Now $\za/{\mathfrak m}\za$ is a standard graded polynomial ring over $K$ in
$n=\mu(\Omega_{A/k})$ variables.
The equality
(\ref{zerointersection}) implies that ${\cal T}$ is mapped
isomorphically onto its image in this polynomial ring.
 Let $h(t)=\dim _K{\cal  T}_t$ and $r=\min\{t\geq 0\,|\, {\cal T}_t\neq
0\}$). It follows that $h(t)$
is at least the number of monomials of degree $t-r$, in other words,
\[h(t)\geq \left( \begin{array}{c}
                     t-r+n-1\\
                       n-1
\end{array} \right).\]

Now consider a presentation $A\simeq R/I$ as given by assumption,
where $(R,{\mathfrak n})$ is a regular local ring essentially of finite type over $k$ and
$I\subset{\mathfrak n}^2$. Set $m=\mu
(I+{\mathfrak n}^3/{\mathfrak n}^3)$ and
$(A',{\mathfrak m}')=(A/{\mathfrak m}^2,{\mathfrak m}/{\mathfrak m}^2)$.
The usual $A$-free presentation of
$\Omega_{A/k}$ by means of a Jacobian matrix induces a presentation
\[A'\,^m\lar A'\,^n\lar \Omega_{A/k}\otimes_A A' \lar 0,\]
yielding, for every $t\geq 1$, an exact sequence
\[A'\,^m\otimes_{A'} {\cal S}_{t-1}(A'\,^n)\lar {\cal S}_t(A'\,^n)\lar
{\cal S}_t(\Omega_{A/k})\otimes_A A'\lar 0.\]

The equality (\ref{zerointersection}) also says that, for every $t\geq 0$, the
graded piece ${\cal T}_t$ of ${\cal T}$ is a $K$-vector space of dimension $h(t)$
and a direct
summand of ${\cal S}_t(\Omega_{A/k})$ as an $A$-module. Hence ${\cal S}_t(\Omega_{A/k})\otimes_A A'$
too admits $K^{\oplus h(t)}$ as a direct summand. Therefore
${\mathfrak m}'\,^{\oplus h(t)}$ is a direct summand of the image of
$A'\,^m\otimes_{A'} {\cal S}_{t-1}(A'\,^n)$ in $
{\cal S}_t(A'\,^n)$, which implies that
\[m \left( \begin{array}{c}
                     t-1+n-1\\
                       n-1
\end{array} \right) \geq
\mu({\mathfrak m}')h(t)\geq \mu({\mathfrak m}') \left(
\begin{array}{c}
                     t-r+n-1\\
                       n-1
\end{array} \right),\]
for every $t$. Observing that $\dim R=\mu({\mathfrak
m})=\mu({\mathfrak m}')$, this inequality contradicts the
assumption $m\leq \dim R -1$. \QED

\begin{Corollary}\label{quadrics}
Let $(R,{\mathfrak n})$ be a regular
local ring essentially of finite type over a perfect field $k$, and let
$I\subset R$ be an ideal such that $A=R/I$ is reduced. Assume that
one of the following conditions holds\/$:$

\noindent {\rm (a)} $\mu(I_{{\mathfrak p}})\leq \dim R_{{\mathfrak p}}-1\ $ for every
non-minimal  ${\mathfrak p}\in V(I)${\rm ;}
\newline{\rm (b)} $I\subset {\mathfrak n}^3$ and $A$ is an isolated
singularity.

Then $\za$ is reduced {\rm (}if and{\rm )} only if it is $A$-torsionfree.
\end{Corollary}

\begin{Remark}\label{superquadrics} \rm The full force of
Theorem~\ref{redtorsionfree} takes form in the case of a homogeneous
ideal $I\subset k[\XX]=k[X_1,\ldots,X_n]$ such that $A=k[\XX]/I$ is
reduced and regular on the punctured spectrum. If $\za$ is reduced,
but {\it not\/} torsionfree, then $\mbox{\rm Proj}\,
(k[\XX]/I)\subset \pp_k^{n-1}$ lies on the intersection of $n$
independent quadrics because in this situation $\dim_k[I]\,_2 
\geq n$ by Theorem~\ref{redtorsionfree}. There is a converse when
$\dim_k[I]\,_2\leq 2 \, \hht \,I$, namely, if $\za$ is torsionfree
then $n\geq 2 \, \hht \,I+1$ (cf. the beginning of Section~\ref{RAA}),
hence $\mbox{\rm Proj}\, (k[\XX]/I)\subset \pp_k^{n-1}$ cannot lie
on the intersection of $n$ independent quadrics.
\end{Remark}

The next two examples illustrate the above results. The first shows that
Theorem 2.1 is sharp, i.e., that the
assumption on the numbers of generators cannot be relaxed.

\begin{Example}\rm
Consider the homogeneous coordinate ring $A$ of the
Veronese surface in $\pp_k^5$. Here $A=k[\XX]/I_2(\XX)$ where $\XX$ is a
symmetric $3\times 3$ matrix of indeterminates over $k$
(assumed to be of characteristic $\neq 2$). Notice that $A$ is a Cohen-Macaulay
domain with an isolated singularity, $I_2(X)$ is generated by quadrics,
and $\mu(I_2(X))= 6 = \dim k[X]$.
It can be seen
with the aid of {\it Macaulay\/} that $\za$ is reduced, but not
$A$-torsionfree. In fact, $\za$ has exactly two minimal primes,
${\cal T}$ and $(\XX)\za$, the first of which defines $\ra$.

Another computation with {\it Macaulay\/} shows that $\za$ is Cohen--Macaulay. This
example incidentally answers a  question posed in \cite[Remark after (2)]{HSV1}
concerning the existence of a Cohen--Macaulay generic complete intersection
 ideal ${\bf D}\subset  S$ for which the symmetric algebra ${\cal S}({\bf D}/{\bf D}^2)$ is
Cohen--Macaulay, but the symmetric algebra ${\cal S}({\bf D})$ is not.
Mark Johnson obtained examples of the latter behavior even in
a polynomial ring.
\end{Example}

\begin{Example}\rm Let $A=k[\XX]/I_2(\XX)$ where $(\XX)$ is the
$r$-catalecticant  matrix
$$\left(\begin{array}{cccc}
X_1 & X_2 & X_3 & X_4\\
X_{r+1} & X_{r+2} & X_{r+3} & X_{r+4}
\end{array}
\right)
$$
with $1\leq r\leq 4$. Note that for $r=4$ we obtain the $2\times 4$
generic  matrix and for $r=1$, the usual $2\times 4$ generic Hankel
matrix. Also $A$ is a Cohen--Macaulay (domain) in all cases, being a
specialization of the generic situation. Clearly $A$ is an isolated
singularity. For values $1\leq r\leq 2$ the Zariski algebra $\za$ is
neither Cohen--Macaulay nor torsionfree. Geometrically, we see the
reason for $\za\neq \ra$ since in this range the tangential variety
$\mbox{\rm Proj}\, (k\otimes_A \ra)  \subset \pp_k^{r+3}$
to $\mbox{\rm Proj}\, (k[\XX]/I_2(\XX))\subset \pp_k^{r+3}$ is a proper
subvariety of the ambient space. For values $3\leq r\leq 4$ the
Zariski algebra $\za$ is Cohen--Macaulay (computer calculation),
hence torsionfree because $r+4\geq 2\,\hht \,I_2(\XX)+1=7$ (see
\cite[3.3]{sv}).
\end{Example}

\subsection{Analogues of Berger's conjecture}

When $\dim A=1$ one can make
Theorem~\ref{redtorsionfree} more precise.
Here it is natural to work in the more general setting of complete
$k$-algebras. Instead of the module of K\"ahler differentials
we use the universally finite module of differentials, defined in an analogous way as
$\dd/\dd^2$ where $\dd$ is the kernel of the multiplication map $A\widehat{\otimes}_k A\lar A$.
Our result is reminiscent of Berger's well-known conjecture asserting that the
 module of differentials of a reduced curve singularity over a perfect field
 cannot be torsionfree (see \cite{Berger}). It is clear at least that the entire
symmetric algebra $\za$ cannot be $A$-torsionfree, as otherwise the module of differentials
would satisfy condition $(F_1)$,
which translates into the inequality $\edim A \leq 2 \dim A -1=1$ (see the discussion
at the beginning of Section~\ref{RAA}).

\begin{Theorem}\label{Bergerlike} Let $k$ be an algebraically
closed field and assume that $A =k[[x_1, \ldots, x_n]]$ is a
one--dimensional domain.  Then $\za$ is reduced {\rm (}if and{\rm )}
only if $A$ is regular.
\end{Theorem}

\demo Assume that $\za$ is reduced. To argue by way of
contradiction assume that $n = \edim A \geq 2$.  Letting $k[[t]]$
denote the integral closure of $A$,  we may arrange so that in the
$t$--adic valuation one has
\[ v(x_1) < v(x_2) < v(x_3) \leq \cdots \leq v(x_n) \]
and $v(x_2)/v(x_1)$ is not an integer.  Write $A = k[[X_1,
\ldots,X_n]]/(f_1, \ldots, f_m)$, with $f_i \in (X_1, \ldots,
X_n)^2$, $T_j=dX_j$ and $t_j=dx_j$. Notice that
\[ \za \simeq A[t_1, \ldots, t_n] :\,= A[T_1, \ldots, T_n]/(\sum_{j=1}^n
\frac{ \partial f_i}{\partial x_j} T_j \ | \ 1 \leq i \leq m).  \]
Here $ \partial f_i/\partial x_j$ denotes the image of
 $\partial f_i/\partial X_j$ in $A$.
Clearly these elements are contained in the maximal
ideal $\fm$ of $A$.

Next, let $J$ be an $A$--ideal isomorphic to $\Omega_{A/k}$ modulo
its torsion and let ${\cal T}$ be the kernel of the natural map from
$\za$ onto the Rees algebra ${\cal R}_A(J)$. By the reducedness of
$\za$, one has ${\cal T} \cap {\fm}\za= 0$ as in the previous
section (see (\ref{zerointersection})). In particular ${\rm
Supp}({\fm}\za) \cap V({\fm}\za) \subset V({\cal T} + {\fm}\za)$ as subsets
of $\spec (\za)$. On the other hand,
\[\dim \za/({\cal T} + {\fm}\za) = \dim {\cal R}_A(J)/\fm {\cal
R}_A (J) = 1, \] and therefore $({\fm}\za, t_3, \ldots,
t_n) \not \in {\rm Supp}({\fm}\za)$.

Thus $({\fm}\,{\za})_{({\fm}\,\za,t_3, \ldots, t_n)} = 0$.
From the presentation
\[\za/(t_3, \ldots, t_n) = A[T_1, T_2]/L \]
with
\[ L = \left(\frac{\partial f_i}{\partial x_1}T_1 + \frac{\partial
f_i}{\partial x_2}T_2 \ | \ 1 \leq i \leq m\right), \]
we see that
\[\fm A[T_1,T_2]_{\fm A[T_1, T_2]} = LA[T_1, T_2]_{{\bf
m}A[T_1, T_2]}.  \]
In  particular, there exists a polynomial $g(T_1, T_2) \in A[T_1,
T_2] \setminus \fm A[T_1, T_2] $ with
\[ g(T_1, T_2)\cdot \fm  \subset  \left(\frac{\partial f_i}{\partial
x_1}T_1 + \frac{\partial f_i}{\partial x_2}T_2 \ | \ 1 \leq i \leq
m\right).  \] Comparing coefficients, we conclude that
\[ \fm  = \left(\frac{\partial f_i}{\partial x_1}, \frac{\partial
f_i}{\partial x_2} \ | \ 1 \leq i \leq m\right).  \]

Since $v(x_1)$ is the smallest positive element in the value
semigroup of $A$, there exists an $i$ such that $v(x_1) =
v(\partial f_i/\partial x_1)$, or  $v(x_1) = v(\partial f_i/\partial
x_2)$. Write $f=f_i$. Assuming the first,  we have $f = aX_1^2 + g$, with $a \in k$
and $g \in (X_1^3) + (X_1, \ldots, X_n)(X_2, \ldots, X_n)$.
Obviously $v(g(x_1, \ldots, x_n)) > 2v(x_1)$ and hence $v(\partial
g/\partial x_1)
>v(x_1)$. Therefore $v(\partial f/\partial x_1) = v(x_1)$ implies $a
\neq 0$, while $f(x_1, \ldots, x_n) = 0$ implies $a = 0$. This is a
contradiction.

Next assume $v(\partial f/\partial x_2) = v(x_1)$ and write
\[f= \sum_{i = 2}^{\infty}a_i{X_1}^i + bX_1X_2 + h, \]
where $a_i \in k, \ b \in k$ and $h \in (X_1^2, X_2)X_2 + (X_1,
\ldots, X_n)(X_3, \ldots, X_n)$. Here $v(h(x_1, \ldots, x_n)) >
v(x_1) + v(x_2)$ and hence $v(\partial h/\partial x_2) > v(x_1)$.
Therefore $v(\partial f/\partial x_2) = v(x_1)$ implies $b \neq 0$.
But then $f(x_1, \ldots, x_n) = 0$ implies
\[v(x_1) + v(x_2) = v(bx_1x_2) =
v(\sum_{i=2}^{\infty} a_ix_1^i) = lv(x_1) \]
for some integer $l$.  But this
contradicts our assumption that $v(x_2)/v(x_1)$ is not an
integer.  \QED

\begin{Remark}\rm The above result means that the first jet scheme of an algebroid
curve singularity over an algebraically closed field cannot be reduced.
\end{Remark}

The proof of Theorem~\ref{Bergerlike} also gives the following stronger result:

\begin{Theorem}\label{Bergerlike2} Let $k$ be a perfect field, let $A = k[[x_1,
\ldots, x_n]]$ be a one--dimensional ring and assume that  $\za$ is reduced. Then
$A/{\mathfrak p}$ is regular  for every minimal prime ideal  ${\mathfrak p}$ of $A$.
If in addition ${\rm char}(k) = 0$ and $A$ is quasi--homogeneous {\rm (}i.e.,
the completion of a positively graded $k$--algebra{\rm )}, then $A$ is regular.
\end{Theorem}

\demo We use the notation of the previous proof.  Of course $A$
is reduced and we may assume that $k$ is algebraically closed.
Suppose there exists a minimal prime ${\mathfrak p}$ of $A$ such that
$A/{\mathfrak p}$ is not
regular.
Then choose a minimal generating set $x_1, \ldots, x_n$ of $\fm$
such that \[v(x_1+{\mathfrak p}) < v(x_2 + {\mathfrak p}) < v(x_3+ {\mathfrak p})
 \leq \cdots\leq v(x_n+{\mathfrak p})
\leq \infty\] and the ratio of the first two integers is not an
integer (here $v$ is the $t$--adic valuation on  the
integral closure $k[[t]]$ of $A/{\mathfrak p}$).
 The rest of the proof proceeds as before replacing $\za/(t_3, \ldots, t_n)$
by $\za/({\mathfrak p} S, t_3, \ldots, t_n)$.

Now assume that ${\rm char} (k) = 0$ and $A$ is quasi--homogeneous.
  Let $\tau$ denote the kernel of the
Euler map $\Omega_{A/k}\surjects \fm$.  Since $\tau$ is the torsion of $\Omega_{A/k}$,
by the reducedness of $\za$
we have $\tau \bigcap \fm \Omega_{A/k} = 0$ as in the proof of
Theorem~\ref{redtorsionfree},
hence $\tau$ is a direct summand
of $\Omega_{A/k}$. It follows that $\tau = 0$ since $\mu(\Omega_{A/k})
= {\rm edim \ } R = \mu(\fm)$. But  then $A$ is regular because in
the quasi--homogeneous case the Berger conjecture is true
(see \cite[4.4]{SS}). \QED

\bigskip

We shall now consider  another version of Berger's conjecture.
Throughout $\tau_A(E)$ will denote the $A$-torsion of an $A$-module $E$.
If $(A,{\fm})$ is a Noetherian local ring, we write $\ecodim A=\mu(\fm)-\dim A$.
Let $A\simeq R/I$ be a Cohen--Macaulay normal ring, where $R$ is a regular local ring
essentially of finite type over a perfect field $k$ and $I$ is an ideal of height $2$
(i.e.,  $\ecodim A \leq 2$).
Here the first Koszul homology ${\rm H}_1(I)$ of a generating set of $I$
is Cohen--Macaulay (see \cite[2.1(a)]{AvHe1}), in particular the natural complexes
\begin{equation}\label{conormalsequence}
 0\lar {\rm H}_1(I)\lar A^m{\lar} I/I^2\lar 0
 \end{equation}
 and
 \begin{equation}\label{omegasequence}
 0\lar  I/I^2\lar  A^n\lar \Omega_{A/k}\lar 0
 \end{equation}
 are exact.

\begin{Theorem}\label{Bergerlite}
Let $(A,\fm)$ be a normal local Cohen--Macaulay ring essentially of finite type
over a perfect field $k$ satisfying the following conditions{\rm :}
\begin{enumerate}
\item[{\rm (a)}] $\ecodim A=2${\rm ;}
\item[{\rm (b)}] $A$ is an almost complete intersection locally in codimension $2$.
\end{enumerate}
If $\,\Omega_{A/k}/\tau_A(\Omega_{A/k})$ is a Cohen--Macaulay module
then $A$ is regular.
\end{Theorem}
\demo Set $\Omega=\Omega_{A/k}$ and let $\omega=\omega_A$ denote the
canonical module of $A$.  Write $A=R/I$, where $R$ is a regular local
ring of essentially of finite type and of transcendence degree $n$ over $k$
and $I$ is a perfect ideal of height $2$.
It suffices to prove that $I$ is a complete intersection. For in this case $\Omega$
is torsionfree of finite projective dimension -- thus $\Omega=\Omega/\tau_A(\Omega)$
is a maximal Cohen--Macaulay module of finite projective dimension, hence
necessarily free.

We first consider the case $\dim A=2$.
Supposing $I$ is not a complete intersection we have that this ideal is
 minimally generated by $3$ elements.
 Consider a minimal presentation
 \begin{equation}\label{presentation}
 0\lar R^2\stackrel{\phi}{\lar} R^3\lar I\lar 0.
 \end{equation}
 From (\ref{conormalsequence}) and (\ref{omegasequence}) we obtain an exact sequence
 \begin{equation}\label{canonical}
 0\lar \omega\simeq {\rm H}_1(I)\lar A^3\stackrel{\Theta}{\lar} A^n\lar \Omega\lar 0,
 \end{equation}
 with $\Theta$ the transposed Jacobian matrix of
 the $3$ generators of $I$.
 Moreover, since $\Omega/\tau_A(\Omega)$ is Cohen--Macaulay and $\tau_A(\Omega)$ has grade
 at least $2$,
${\rm Ext}^1_A(\Omega,\omega)=0$.
 Thus applying $-^{\vee}={\rm Hom}_A(-,\omega)$ to (\ref{canonical}) one can see that
 there is an exact sequence
$$
\omega^n\stackrel{\Theta^{\vee}}{\lar} \omega^3\lar \omega^{\vee}\simeq A
 \lar A/I_1(\phi)A\lar 0.
$$
Hence we obtain the exact sequence
\begin{equation}\label{canonicalbis}
\omega^n\stackrel{\Theta^{\vee}}{\lar} \omega^3 \lar I_1(\phi)/I\lar 0.
\end{equation}

As $I$ is  a perfect ideal of height $2$, but not a complete intersection, one has
$I\subset I_1(\phi)^2$.
In particular, $I_1(\Theta)\subset I_1(\phi)$.
Thus tensoring (\ref{canonicalbis}) with $R/I_1(\phi)$ gives
\[ (\omega/I_1(\varphi)\omega)^3 \simeq I_1(\varphi)/I \otimes_R
R/I_1(\varphi)\simeq I_1(\varphi)/I_1(\varphi)^2. \]
On the other hand, dualizing (\ref{presentation}) into $R$ yields
$\omega\otimes_R R/I_1(\varphi)\simeq (R/I_1(\varphi))^2$. It
follows that
\[ I_1(\varphi)/I_1(\varphi)^2\simeq
 (\omega/I_1(\varphi)\omega)^3\simeq (R/I_1(\varphi))^6.\]
 Since $R$ is regular, \cite[1.1]{Va0} implies that $I_1(\varphi)$ is
generated by a regular sequence of length $6$. But this is impossible
in $R$, a ring of dimension four.

We now consider the case of arbitrary dimension.
By the above, $A$ is regular locally in codimension $2$.
Therefore (\ref{conormalsequence}) and (\ref{omegasequence}) imply
that $\Omega$ is torsionfree, hence Cohen--Macaulay by
assumption.
But then also $I/I^2$ is Cohen--Macaulay by (\ref{omegasequence}), which forces
$I$ to be a complete intersection (see \cite[2.4]{He78}).
\qed

\subsection{Normal Zariski tangent algebras}

In this part we give sufficient conditions for the normality of $\za$ to imply
the reflexivity of it graded components. The conditions are stated in terms of
$\mu(I+{\mathfrak n}^3/{\mathfrak n}^3)$, where $I$ is the defining ideal
of $A$ in $(R,\fn)$.
We will write $\ell(\cdot)$ for the analytic spread of an 
$R$-ideal (see also the beginning of Section 3).

\begin{Proposition}\label{normality_quadrics}
Let $(R,{\mathfrak n})$ be  a regular local ring essentially of finite type over a
perfect field $k$ and let $I\subset {\mathfrak n}^2$ be an ideal such that $A=R/I$ is
reduced.
Write $g$ for the height of $I$ and $\fm $ for the maximal ideal of $A$.
\begin{enumerate}
\item[{\rm (a)}] If $\za$ is equidimensional and $(\za)_{(\fm)}$ is regular, then
$\mu(I+{\mathfrak n}^3/{\mathfrak n}^3)\geq 2g$.
\item[{\rm (b)}] Assume that ${\rm char} (k)=0$ and 
$R$ is the local ring of the polynomial ring $k[X_1,\ldots,X_n]$ at its
homogeneous  maximal ideal. Let $I_2$ denote the
ideal of $R$
generated by $2$-forms in $X_1,\ldots,X_n$ such that $I+{\mathfrak n}^3=I_2+{\mathfrak n}^3$.
If $\za$ is equidimensional and $(\za)_{(\fm)}$ is regular, then $\ell(I_2)=2g$.
\end{enumerate}
\end{Proposition}
\demo
Knowingly,  ${\mathbb S}_{R/k}$ is a polynomial ring over $R$. More precisely,
${\mathbb S}_{R/k}\simeq R[\TT]$, where $\TT=T_1,\ldots,T_n$ may be chosen to be the
differentials
of a transcendence basis of $R$ over $k$ in part (a) and the differentials of
$X_1,\ldots,X_n$ in part (b).
Set $K=R/{\mathfrak n}$, so that the residue field of $R(\TT)$ is $K(\TT)$.
Write $d:=\trdeg _k(R/\fn)=\trdeg _k(A/\fm)=\trdeg _k(K)$.
There is a presentation
$$(\za)_{(\fm)}\simeq R(\TT)/{\cal J} .$$
Since $\za$ is equidimensional, we have
\begin{eqnarray*}
\hht \, {\cal J}&=&\dim R(\TT)-\hht \,\fm\za=d+n-(\dim \za-n)\\
&=&d+n-((d+2n-2g)-n))=2g,
\end{eqnarray*}
and since $(\za)_{(\fm)}$ is regular, we have
$$\dim _{K(\TT)} ({\cal J}+{\mathfrak n}^2/{\mathfrak n}^2)=2g.$$

Now let $f_1,\ldots,f_m\in R$ be so chosen that
they generate $I$ modulo ${\mathfrak n}^3$ for part (a), and for part (b) they be
homogeneous quadrics in
$k[X_1,\ldots,X_n]$ generating the ideal $I_2$.
In either case, one has
$${\cal J}+{\mathfrak n}^2/{\mathfrak n}^2=
(df_i\,|\, 1\leq i\leq m)+{\mathfrak n}^2/{\mathfrak n}^2.$$
It follows that $m\geq 2g$, thus proving (a).

To prove (b), let $f=f(\XX)$ denote any of the chosen forms
$f_1,\ldots,f_m$ in $k[\XX]=k[X_1,\ldots,X_n]$. Then
$$\XX\cdot \phi\cdot \TT^t=\TT\cdot \phi\cdot \XX^t,$$
where $\phi$ stands for the Hessian matrix of $f$. From
 this one sees that
$$\sum_{j=1}^n \frac{\partial f(\XX)}{\partial X_j}T_j=
\sum_{j=1}^n \frac{\partial f(\TT)}{\partial T_j}X_j,$$
hence
$$df=\sum_{j=1}^n \frac{\partial f(\TT)}{\partial T_j}X_j.$$
We deduce that
\begin{eqnarray*} \dim _{k(\TT)}(df_i\,|\, 1\leq i\leq m)+{\mathfrak n}^2/{\mathfrak n}^2&=&
\rank \left(\frac{\partial f_i(\TT)}{\partial T_j}\right)_{i,j}
=\rank \left(\frac{\partial f_i(\XX)}{\partial X_j}\right)_{i,j}\\
&=& n-\rank \Omega_{R/k[f_1,\ldots,f_m]}=n-\trdeg _{k[f_1,\ldots,f_m]}\,R\\
&=& \trdeg_k \,k[f_1,\ldots,f_m]=\dim k[f_1,\ldots,f_m]
\end{eqnarray*}
(see also \cite[1.1]{SimisDiff}).
On the other hand, since $f_1,\ldots,f_m$ are forms of the same degree generating
the ideal $I_2$,
we have $k[f_1,\ldots,f_m]\simeq {\cal R}(I_2)\otimes _R k$, hence $\dim k[f_1,\ldots,f_m]
=\ell(I_2)$.
Summing up, we have proved that
$$2g=\dim _{k(\TT)} ({\cal J}+{\mathfrak n}^2/{\mathfrak n}^2)=\ell(I_2).$$
\qed

\begin{Theorem}\label{normality_reflexiveness}
Let $A$ be a normal domain essentially of finite type over a perfect field $k$.
Assume that one of the following conditions holds $:$
\begin{enumerate}
\item[{\rm (a)}] For every non-regular prime ${\mathfrak p}\in
\spec (A)$, $A_{{\mathfrak p}}\simeq R/I$
where $(R,{\mathfrak n})$ is a regular local ring essentially of finite type
over $k$ and $I\subset
{\mathfrak n}^2$ is an $R$-ideal satisfying
$\mu(I+{\mathfrak n}^3/{\mathfrak n}^3)\leq 2 \, \hht \, I-1$;
\item[{\rm (b)}]  $\mbox{\rm (char}(k)=0${\rm )} $(R,{\mathfrak n})$ is the
local ring of $k[X_1,\ldots,X_n]$ at $(X_1,\ldots,X_n)$,
$A=R/I$, and $\ell(I_2)\neq 2 \, \hht \, I$, where $I_2$ denotes the ideal of $R$
generated by $2$-forms in $X_1,\ldots,X_n$ such that
$I+{\mathfrak n}^3=I_2+{\mathfrak n}^3$.
\end{enumerate}
Then $\za$ is normal {\rm (}if and{\rm )} only if ${\cal
S}_i(\Omega_{A/k})$ is a reflexive module for every $i\geq 1$.
\end{Theorem}
\demo We prove that the natural inclusion $\za\subset
\bigoplus_{i\geq 0}{\cal S}_i(\Omega_{A/k})^{**}$ is an equality
assuming that $\za$ is normal. It suffices to show this after
localizing at the codimension one primes of $\za$. Thus, let
${\mathfrak q}\in \spec (\za)$ be a prime of height one. We are done
once we prove that $A_{{\mathfrak q}\cap A}$ is regular. Now,
replacing $A$ by $A_{{\mathfrak q}\cap A}$, we may assume that
$(A,\fm)$ is local and $\fm\subset {\mathfrak q}$. By assumption,
$(\za)_{\mathfrak q}$ is regular, hence so is its further
localization $(\za)_{(\fm)}$. According to
Proposition~\ref{normality_quadrics} this contradicts our
assumptions, unless $A$ is regular. \qed

\section{The Rees algebra of the module of differentials}\label{RAA}

We will draw upon \cite{ram1} for terminology and basic notions
about Rees algebras of modules.  Let $E$ be a finite module over a Noetherian ring
$A$ and assume that $E$ is generically free (i.e., free locally at every associated
prime of $A$). In this setting
the {\it Rees algebra\/} ${\cal R}_A(E)$ of $E$ is defined as the symmetric algebra
of $E$ modulo $A$-torsion,  ${\cal R}_A(E)={\cal S}_A(E)/\tau_A({\cal S}_A(E))$.
As is known,
this definition retrieves the usual notion of the Rees algebra of an
ideal containing a regular element.
If in addition $E$ is torsionfree, then $E$ can be embedded into a
finite free module $F$ and the Rees algebra of $E$
can be identified with the image of ${\cal S}_A(E)$ in ${\cal S}_A(F)$.
In case ${\cal S}_A(E)={\cal R}_A(E)$ we say that $E$ is a module of {\it linear type}.

We recall the following additional notions. Let $E$ be a finite
module over a Noetherian ring $A$. Suppose that $E$ has a rank $r$
(i.e., $E$ is free of rank $r$ locally at every associated prime of
$A$). Given an integer $t\geq 0$, we say that $E$ satisfies
condition $(F_t)$ if $\mu(E_{\mathfrak{p}})\leq r+ \dim A_{\mathfrak
p} -t$ for every $\mathfrak{p}\in \spec(A)$ such that
$E_\mathfrak{p}$ is not free. In terms of Fitting ideals this
condition is equivalent to the condition that $\height \mbox{\rm
Fitt}_i(E)\geq i-r +t+1$ for $ i\geq r$.

By the same token, one can introduce yet another condition  based on
the analytic spread $\ell(E)$ of a module over a local (or graded)
ring $(A,\fm)$, defined to be the Krull dimension of the residue
algebra ${\cal R}_A(E)/\fm {\cal R}_A(E)$ (see \cite{ram1}). We will
say that a finite module $E$ over a Noetherian ring $A$ satisfies
condition $(L_t)$ to mean that $\ell(E_{\mathfrak p})\leq \dim
A_{\mathfrak p}+r-t$ for every ${\mathfrak p}\in \Spec(A)$ with
$\dim A_{\mathfrak p}\geq t$. Roughly, this condition plays a
similar role for the Rees algebra as $(F_t)$ plays for the symmetric
algebra. We will only use these conditions in the range $0\leq t\leq
2$.

\medskip

We will mainly focus on the case where $E=\Omega_{A/k}$. Let $A$ be
a reduced algebra essentially of finite type over a perfect field
$k$ -- hence $\Omega_{A/k}$ is generically free. We denote the Rees
algebra of $\Omega_{A/k}$ by $\ra$. If in addition $(A,\fm)$ is
equidimensional then $\Omega_{A/k}$ has a rank and this rank equals
 $\dim A+\mbox{\rm trdeg}_k(A/\fm)$.
In this case some of the above Fitting conditions can be expressed in terms of local
embedding dimensions.
Namely,
$\Omega_{A/k}$ satisfies $(F_t)$ if and only if $\edim A_{\mathfrak p}
 \leq 2\dim A_{\mathfrak p}-t$ for every
non-regular prime ${\mathfrak{p}}\in {\rm Spec}(A)$ (see \cite[the proof of 2.3]{star}).
 The result of the next subsection will be stated in these terms.

\subsection{Cohen--Macaulayness of $\mathbb{R}_{A/k}$}
\label{Reescm}

Let $A$ be a local Cohen--Macaulay ring and let $E$ be a finite
$A$-module of projective dimension one. If $E$ satisfies condition
$(F_1)$ then the symmetric algebra ${\cal S}_A(E)$ is a
Cohen--Macaulay torsionfree $A$-algebra (see \cite[Proposition
4]{Avramov}, \cite[1.1]{Huneke}, \cite[3.4]{sv}) -- in particular,
the Rees algebra ${\cal R}_A(E)$ is Cohen--Macaulay. The question as
to whether, conversely, ${\cal R}_A(E)$ being Cohen--Macaulay
implies the condition $(F_1)$ fails in general even if $E$ satisfies
$(F_0)$ (see \cite[4.7]{ram1}). The theorem below will show that
this converse holds in the case of the module of differentials of a
complete intersection over a field of characteristic zero. This
shows that Cohen--Macaulayness is a rather restrictive property for
${\mathbb R}_{A/k}$.

\begin{Theorem}\label{cmversusF1}
Let $k$ be a field of characteristic zero and let $A$ be a
$k$-algebra essentially of finite type which is locally a complete
intersection. Assume the following$\,${\rm :}
\begin{enumerate}
\item[{\rm (i)}] $\edim A_{\mathfrak{p}}\leq 2\dim A_{\mathfrak{p}} $
for every prime ${\mathfrak{p}}\in {\rm Spec}(A)\,${\rm ;}
\item[{\rm (ii)}] $\mathbb{R}_{A/k}$ is Cohen--Macaulay.
\end{enumerate}
Then $\edim A_{\mathfrak{p}}\leq 2\dim A_{\mathfrak{p}}-1$ for every
non-minimal prime ${\mathfrak{p}}\in {\rm Spec}(A)$.
\end{Theorem}

\demo Arguing by way of contradiction, let ${\mathfrak{p}}\in {\rm
Spec}(A)$ be minimal such that $\edim A_{\mathfrak{p}}= 2\dim
A_{\mathfrak{p}} \geq 2$. By localizing at this prime, we may assume that
$(A, {\mathfrak m})$ is local and that $\edim A=2\dim A$.

Next we reduce the argument to the case where the residue field
$A/{\mathfrak m}$ is algebraic over $k$. To do this, let $r$ be the
transcendence degree of $A/\fm$ over $k$ and suppose $r\geq 1$.
Write $A=k[x_1,\ldots,x_m]_{\mathfrak p}$ and pick $r$ general
$k$-linear combinations $y_1,\ldots,y_r$ of $x_1,\ldots,x_m$ such
that, in particular, their residues yield a transcendence basis of
$A/{\mathfrak m}$ over $k$. Note that $K:=k(y_1,\ldots,y_r)$ is a
subfield of $A$. Furthermore,
$\Omega_{A/K}=\Omega_{A/k}/(Ady_1+\cdots+Ady_r)$. Hence, by the
general choice of $y_1,\ldots,y_r$, the Rees algebra
$${\cal R}(\Omega_{A/K})\simeq {\cal R}(\Omega_{A/k})/(dy_1,\ldots,dy_r)$$ is
again Cohen--Macaulay (see \cite[2.2(f)]{CPU}). Thus, replacing $k$ by $K$, we may assume
henceforth that $A/{\mathfrak m}$ is algebraic over $k$.

Set $d=\dim A$ and induct on $d$. If $d=1$, the Cohen--Macaulayness
of $\mathbb{R}_{A/k}$ implies that $\Omega_{A/k}$ modulo torsion is
free (see \cite[4.3]{ram1}), which would imply that $A$ is regular
(see \cite{Lipman}), thus contradicting the equality $\edim A=2$.

Thus we may suppose that $d\geq 2$. Write $A=R/I$, where $(R,{\mathfrak n})$ is a
regular local ring
essentially of finite type over $k$ with $n:=\dim R=\edim A=2d$. Recall
that $n=\mbox{\rm trdeg}_kR$, hence $\Omega_{R/k}=
RdX_1\oplus \cdots\oplus RdX_n$.
 By assumption, $I$ is generated by an
$R$-regular sequence $f_1, \ldots, f_{d}$. Denoting the images of
$\frac{\partial f_j}{\partial X_i}$ in $A$ by $\frac{\partial
f_j}{\partial x_i}$ we consider the $n\times d$ matrices
\[ \Theta = (\frac{\partial f_j}{\partial X_i}), \qquad \theta =
(\frac{\partial f_j}{\partial x_i}). \] Now $\theta$ presents the
$A$--module $\Omega_{A/k}$, which has projective dimension one, rank
$d$, and satisfies condition $(F_1)$ locally in codimension $d-1\geq
1$ by the inductive hypothesis. Since $\mathbb{R}_{A/k}$ is
Cohen--Macaulay, we may apply \cite[4.6(a) and 4.7]{ram1} with
$s=2d-1$ to show that after a linear change of variables,
\[ I_1(\theta)  = \left(\frac{\partial f_1}{\partial x_n}, \ldots,
\frac{\partial f_{d}}{\partial x_n}\right). \]
Thus
\[ I_1(\Theta)  \subset \left(\frac{\partial f_1}{\partial X_n}, \ldots,
\frac{\partial f_{d}}{\partial X_n}, \,I\right). \]

Set $J=I_1(\Theta)$ and $J_0=(\frac{\partial f_1}{\partial X_n}, \ldots,
\frac{\partial f_{d}}{\partial X_n})$, and notice that $J+I=J_0+I$ by the above.
Since $R/{\mathfrak n}$ is algebraic over
$k$ and char$(k)=0$, we have $I\subset \overline{{\mathfrak n}J}$,
where $\,\bar{}\,$ denotes integral closure (see \cite[Exercise 5.1]{Huneke2} in the case
where $R$ is a power series ring over $k$). Therefore
$$ J + I = J_0 +I \subset J_0 +\overline{{\mathfrak n}J}\subset
\overline{J_0+{\mathfrak n}(J+I)}.$$ The valuative criterion of
integrality now yields $J+I\subset \overline{J_0}$. Therefore ${\rm
height} (J+I)=\hht\, J_0\leq d$. Hence $\dim R/(J+I) \geq 2d-d =d \geq 2$.
But this is impossible because $R/(J+I) \simeq A/I_1(\theta)$ has dimension
zero by assumption (i).
\QED

\begin{Remark}\rm The theorem implies the following geometric result:
let $X\subset {\mathbb P}^{2d+1}_k$ denote a $d$-dimensional smooth complete
intersection over a field $k$ of characteristic zero, with homogeneous coordinate ring $A$.
If ${\mathbb R}_{A/k}$ is Cohen--Macaulay then $X$ is degenerate.
\end{Remark}

\subsection{Normality of $\mathbb{R}_{A/k}$}

Let $A$ be a normal domain and let $E$ be a finitely generated
$A$-module of rank $r$. Throughout this and the subsequent sections,
we set $E^*:={\rm Hom}_A(E,A)$ -- the $A$-dual module of $E$. We
write $E^i$ for the $i$th graded piece $\Rees_A (E)_i$ of $\Rees_A
(E)$ and call it the $i$th {\it power\/} of $E$. Notice that
$E^i={\cal S}_i(E)/\tau_A({\cal S}_i(E))$. Let
$\overline{\Rees_A (E)} $ denote the integral closure of $\Rees_A
(E)$ in its field of fractions.
Since $E$ has rank $r$, $E^1$ can be embedded into $A^r$ and any such embedding
induces an embedding of $\Rees_A (E)$ into a polynomial ring $A[{\mathbf t}]=A[t_1,\ldots,t_r]$.
Since $A$ is normal, this further induces an embedding of $\overline{\Rees_A (E)}$ into
$A[{\mathbf t}]$ as a graded $A$-subalgebra.
 We denote by $\overline{E^i}\subset A[{\mathbf t}]_i$ the $i$th
 graded piece of $\overline{\Rees_A (E)}$ and call it the {\it $i$th
 normalized power\/} of $E$.
One can see that there are inclusions
\begin{equation}\label{inclusions} E^i\subset \overline{E^i}\subset
(E^i)^{**}\subset A[{\mathbf t}]_i.
\end{equation}
The algebra ${\cal B}_A(E):=\bigoplus_{i\geq 0}(E^i)^{**}$ is a Krull domain,
not necessarily Noetherian. An important feature of this ring is that it
has the same divisor class group as $A$.
We say that $E$ is {\it integrally closed\/} if
$E^1=\overline{E^1}$.
Finally, $E$ is said to be {\it normal\/} if $E^i=\overline{E^i}$
for every $i\geq 1$, or equivalently, if $\Rees_A (E)$ is normal.

\begin{Proposition}\label{normalityandreflexiveness}
Let $A$ be a universally catenary normal domain and let $E$ be a finitely
generated $A$-module. The following conditions are equivalent{\rm :}
\begin{enumerate}
\item[{\rm (i)}] $\Rees_A (E)$ is normal and $E$ satisfies condition
$(L_2)${\rm ;}
\item[{\rm (ii)}] $\Rees_A (E)$ satisfies condition $(S_2)$ of Serre
and $E$ satisfies condition $(L_2)${\rm ;}
\item[{\rm (iii)}] $\Rees_A (E)={\cal B}_A(E)$.
\end{enumerate}
\end{Proposition}
\demo (i) $\Rightarrow$ (ii) This is obvious.

(ii) $\Rightarrow$ (iii) Let ${\mathfrak p}\subset A$ be a prime ideal of height at
least two.
We may assume that $A$ is local with maximal ideal $\fm={\mathfrak p}$, and that $A/\fm$ is infinite. Then $\hht
\,\fm \Rees_A (E)=\dim \Rees_A (E)-\ell(E)=d+r-\ell(E)\geq 2$, the
first equality by the catenarian assumption and the last inequality by $(L_2)$.
Since $\Rees_A (E)$ satisfies condition $(S_2)$, the ideal $\fm$ contains
a $\Rees_A (E)$-regular sequence of two elements. This is also
a regular sequence on each power $E^i$, hence ${\rm depth}(E^i)\geq
2$ for every $i\geq 1$. This implies that all powers of $E$ are
reflexive.

(iii) $\Rightarrow$ (i) By (\ref{inclusions}), we have equality of
powers and normalized powers throughout. This means that $\Rees_A
(E)$ is normal.

Next let ${\mathfrak p}\subset A$ be a prime ideal of height at least
two. We may assume that $A$ is local, with maximal ideal $\fm={\mathfrak p}$.
To prove condition $(L_2)$ is equivalent to show that $\hht\, \fm
\Rees_A (E)\geq 2$, so it suffices to show that $\fm \Rees_A (E)$
has grade at least two.

To see this, let $a\in \fm\setminus \{0\}$. We may assume that $E$
is torsionfree. Write $r:=\rank E$. Since $A$ is normal, $E$ is free
locally in codimension one.  Hence the $A/(a)$-module $E/aE$ has a
rank and this rank is $r$. Furthermore, the $A$-torsion of ${\cal
S}_A(E)$ vanishes locally at the associated primes of $A/(a)$,
therefore it maps to the $A/(a)$-torsion of ${\cal
S}_{A/(a)}(E/aE)$. It follows that there is an induced surjective
map $\Rees_A (E)\surjects \Rees_{A/(a)} (E/aE)$, hence a surjective
map $\Rees_A (E)/a\Rees_A (E)\surjects \Rees_{A/(a)} (E/aE)$ (see
\cite[4.5]{elam99}). Now, for rank reasons the kernel of this map is
an $A/(a)$-torsion module. On the other hand, since $a$ is regular
on every power $E^i$ and these are assumed to be reflexive,
$E^i/aE^i$ is torsionfree for every $i\geq 1$.  It follows that the
map $\Rees_A (E)/a\Rees_A (E)\surjects \Rees_{A/(a)} (E/aE)$ is an
isomorphism.

Since $A$ is normal of dimension at least two, the grade of $\fm$ is
at least two. Therefore the grade of $\fm/(a)$ in $A/(a)$ is at
least one, i.e., the grade of its extension to $\Rees_{A/(a)}
(E/aE)$ is at least one. By the above isomorphism we then get that
the ideal $\fm \Rees_A (E)/a\Rees_A (E)$ has grade at least one,
hence $\fm \Rees_A (E)$ has grade at least two, as required. \qed

\medskip

In the special case of projective dimension one, we can improve the
above result. For this,
we will rely on a basic property of integrally closed modules.
Let $E$ be an integrally closed submodule of $A^r$. The associated primes
of $A^r/E$ are particularly restricted when $E$ is a module of finite
projective dimension. This phenomenon
 was described in the case of ideals by Burch
(see \cite[Corollary 3]{Burch0} and more recently, using a different
approach, by  Goto and
Hayasaka (see \cite[1.1 and 2.3]{GotoH}). The latter technique
extends to modules as pointed out  by J. Hong (see
\cite[4.8]{icm}).

The central notion here is a
condition, known as $\mathfrak{m}$-fullness:
\begin{definition}\label{m-fullness}\rm
Let $(A, \mathfrak{m})$
be a Noetherian local ring with  infinite residue field and let $E\subset A^r$
be a submodule of rank $r$.
One says that $E$ is $\mathfrak{m}$-{\it full\/} if
there exists an element $x\in \mathfrak{m}$ such that
\[ \mathfrak{m}E:_{A^r} x = E.\]
\end{definition}

It can be shown that an integrally closed module $E\subset A^r$ is $\mathfrak{m}$-full
provided the residue field of the local ring $A$ is infinite.
Most properties of this notion, including Theorem~\ref{BurchGoto} below,
arise from the following result (see \cite[8.101]{intclos} for details).

\begin{Proposition}\label{mfullmod}
Let $(A,\mathfrak{m})$ be a Noetherian local ring and $E$  a
submodule of $A^r$
such that $\mathfrak{m}E :_{A^r} x =E$ for some $x \in \mathfrak{m}$. Then
\[ E/xE \simeq (E :_{A^r} \mathfrak{m} ) /E
\oplus (E+xA^r )/xA^r \simeq (E :_{A^r} \mathfrak{m} ) /E
\oplus E/x(E:_{A^r} \mathfrak{m}). \]
\end{Proposition}

\medskip

\begin{Theorem} \label{BurchGoto}
Let $(A,\mathfrak{m})$ be a Noetherian local ring and $E$  a
submodule of $A^r$
such that $\mathfrak{m}E :_{A^r} x =E$ for some regular element
 $x \in \mathfrak{m}$. Suppose that $\mbox{\rm projdim} (E) <
\infty$ and $\mathfrak{m} \in \mbox{\rm Ass} (A^r /E )$.
 Then $A$ is a regular local ring.
\end{Theorem}

We use this result in the proof of the next theorem. Even for a finite
module $E$ of projective dimension one over normal domain $A$, the normality
of the Rees algebra $\Rees_A(E)$ does not necessarily imply that $E$ satisfies $(F_2)$ -- for
an easy example take $E$ to be the homogeneous maximal ideal of a polynomial ring
in two variables over a field. Surprisingly though, this implication does hold
for modules whose non-free locus is contained in the singular locus of the ring:

\begin{Theorem}\label{projdim1andnormality}
Let $A$ be a Cohen--Macaulay normal domain
and let $E$ be a finitely generated
$A$-module of rank $r$ such that$\,${\rm :}
\newline {\rm (a)} $E$ has a projective resolution $0\rar P_1\lar P_0\lar E\rar 0\,${\rm ;}
\newline {\rm (b)} $E_{\mathfrak p}$ is $A_{\mathfrak p}$-free for any prime ${\mathfrak p}\subset A$
such that $A_{\mathfrak p}\,$
is regular.

The following are equivalent$\,${\rm :}
\begin{enumerate}
\item[{\rm (i)}] $\Rees_A (E)$ is normal$\,${\rm ;}
\item[{\rm (ii)}] $\overline{E^i}=E^i$ in the range $1\leq i\leq
\rank P_0-r\,${\rm ;}
\item[{\rm (iii)}] $E$ satisfies condition $(F_2)$, i.e.,
$\height \mbox{\rm Fitt}_i(E)\geq i-r+3$ for $i\geq r\,${\rm ;}
\item[{\rm (iv)}] $\Rees_A (E)= {\cal B}_A(E)$.
\end{enumerate}
\end{Theorem}
\demo (i) $\Rightarrow$ (ii): This is obvious.

(ii) $\Rightarrow$ (iii): Let
${\mathfrak p}\in \spec (A)$ be a  prime minimal with the property
that $(F_2)$ fails.
Localizing at ${\mathfrak p}$, we may assume that $(A,{\mathfrak m})$ is a local ring of dimension $d$,
$\mu(E)\geq d+r-1$, and $E$ satisfies $(F_2)$ locally on the
punctured spectrum of $A$.
Consider a minimal free resolution
$$0\rar A^{s}\lar A^{n} \lar E\rar 0.$$
By assumption $n\geq d+r-1$, hence $s\geq d-1$.
Since $(F_0)$ holds on the punctured spectrum,
it follows that the Weyman complex is a minimal free resolution
of the symmetric power ${\cal S}_{d-1}(E)$ (see \cite[Theorem 1(b)]{weyman}). This complex
has length $d-1$ because $s \geq d-1$. Moreover
$E$ satisfies $(F_1)$ on the punctured spectrum, therefore
${\cal S}_{d-1}(E)$ is torsionfree.
Thus $E^{d-1}\simeq {\cal S}_{d-1}(E)$ has projective dimension $d-1$.
In particular, $\fm\in\ass ({\cal S}_{d-1}(A^r)/E^{d-1})$.
Since $s\leq \rank P_0-r$, we can apply the
main assumption to conclude that $E^{d-1}$ is integrally closed. Therefore
$A$ is regular by Theorem~\ref{BurchGoto}. This is a
contradiction vis-\`a-vis (b) since $E$ is not free.

(iii) $\Rightarrow$ (iv): This is well-known for modules of projective
dimension one over Cohen--Macaulay rings (see, e.g., \cite[Proposition 4]{Avramov}).

(iv) $\Rightarrow$ (i): This follows from
Proposition~\ref{normalityandreflexiveness}. \qed

\bigskip

Theorem~\ref{projdim1andnormality} applies naturally to the case of complete
intersections via the intervention of the Jacobian criterion, which
says that the non-free locus of the module of differentials coincides
with the singular locus of the ring.

\begin{Corollary}\label{CIandnormality}
Let $A$ be a normal complete intersection domain essentially of finite type over a
perfect field $k$. The following are equivalent$\,${\rm :}
\begin{enumerate}
\item[{\rm (i)}] $\mathbb{R}_{A/k}$ is normal$\,${\rm ;}
\item[{\rm (ii)}] $\overline{\Omega_{A/k}^i}=\Omega_{A/k}^i$ in the range
$1\leq i\leq\, {\rm ecodim}\, A$$\,${\rm ;}
\item[{\rm (iii)}] ${\rm edim}\, A_{\mathfrak p}\leq 2\dim A_{\mathfrak p}-2$ for
every non-regular prime ${\mathfrak p}\in \spec (A)$$\,${\rm ;}
\item[{\rm (iv)}] $\mathbb{R}_{A/k}={\mathbb B}_{A/k}$.
\end{enumerate}
\end{Corollary}
\demo Note that condition (iii) translates into $\Omega_{A/k}$
satisfying $(F_2)$ according to the preliminaries before
Section~\ref{Reescm}. The rest follows immediately from
Theorem~\ref{projdim1andnormality}.
\qed

\bigskip

The following result is reminiscent of Theorem~\ref{redtorsionfree}.

\begin{Proposition} \label{compxdual} Let $A=R/I$ be a normal domain, where $R$ is a
regular ring
essentially of finite type over a perfect field $k$.
If $i\geq 1$ is such that ${\cal S}_i(\Omega_{A/k})\simeq \Omega_{A/k}^{i}$ and
$\Omega_{A/k}^i$ is integrally closed,
 then $I\not \subset \mathfrak{p}^3$ for any $\mathfrak{p}\in
\Ass_R((\Omega_{A/k}^i)^{**}/\Omega_{A/k}^i)$.
\end{Proposition}

\demo We replace $R$ by $R_{\mathfrak p}$ to assume that $(R,{\mathfrak n})$
and $(A,\fm)$ are local. We may also suppose that the residue field is infinite.
Notice that $\fm\in
\Ass_A((\Omega_{A/k}^i)^{**}/\Omega_{A/k}^i)$, in particular $\dim A\geq 2$.
Since ${\cal S}_i(\Omega_{A/k})\simeq \Omega_{A/k}^i$, the presentation
$$
I/I^2 \lar A^n= \Omega_{R/k}\otimes_RA \lar \Omega_{A/k}\rar 0
$$
induces an exact sequence
\begin{equation}\label{higher}
I/I^2 \otimes_A S_{i-1}(A^n) \lar S_i(A^n) \lar  \Omega_{A/k}^i\rar
0.
\end{equation}

Let $d=\trdeg_k A$.
As $\Omega_{A/k}^i$ is assumed to be integrally closed, it is $\fm$-full in any embedding
$\Omega_{A/k}^i\subset A^t$, where $t={{i+d-1}\choose {d-1}}$. Thus there is an element $x\in\fm$
satisfying the hypothesis of Proposition~\ref{mfullmod} with $E=\Omega_{A/k}^i$.
Clearly $\fm\in \Ass_A (A^t/\Omega_{A/k}^i)$, hence $\Omega_{A/k}^i:_{A^t}\fm\neq \Omega_{A/k}^i$.
Therefore the proposition implies that $\Omega_{A/k}^i/x\,\Omega_{A/k}^i$ has
$A/\mathfrak{m}$ as a direct summand over $A/(x)$. On the other hand, tensoring (\ref{higher}) with $A/(x)$,
we see that the syzygies of
 $\Omega_{A/k}^i/x\,\Omega_{A/k}^i$ have coefficients in the $A/(x)$-ideal generated by the
 entries of the
 Jacobian matrix of the generators of $I$. As $\fm/(x)\neq 0$,  we cannot have
 $I \subset \mathfrak{n}^3$.
 \QED

\subsection{Algebras of low codimension}

We now study the normality of $\mathbb{R}_{A/k}$ in several
settings in low dimension or low
embedding codimension.
A highlight is the closer relationship between the normality of
$\mathbb{R}_{A/k}$ and the size of the non-singular locus of $A$ than
is typical of more general Rees algebras.

\medskip

Let $A$ be a local Cohen--Macaulay ring essentially of finite type over
a perfect field $k$ and assume $\ecodim A\leq 2$. Write $A = R/I$, where $R$ is a
regular local ring essentially of finite type over $k$ and $I$ an ideal of height $2$.
 The presentation
\[I/I^2 \lar \Omega_{R/k}\otimes_R A\simeq  A^n \lar \Omega_{A/k} \rar 0\]
induces a complex of graded modules over the polynomial ring $B={\cal S}_A(A^n)$,
\begin{eqnarray}
0  \rar (\wedge^2 I/I^2)^{**} \otimes_A B[-2] \lar I/I^2 \otimes_A
B[-1] \lar B \lar   \mathbb{S}_{A/k} \rar 0.
\label{Koszulcodim2}
\end{eqnarray}
We notice that this is the $\mathcal{Z}$-complex of the module $\Omega_{A/k}$
(see \cite[Chapter 3]{alt} for details).
 It can be used to derive properties of $\mathbb{S}_{A/k}$ and of
$\mathbb{R}_{A/k}$:

\begin{Proposition} \label{normcodim2}
Let $A$ be a reduced local Cohen--Macaulay ring
essentially of finite type over a perfect field $k$.
Assume that $\ecodim A\leq 2$ and that $\edim A_{\mathfrak p}\leq 2\dim A_{\mathfrak p}$ for
every prime ${\mathfrak p}\in \spec (A)$.
\begin{enumerate}
\item[{\rm (a)}] The sequence {\rm (\ref{Koszulcodim2})} is  exact and $\mathbb{S}_{A/k}$
is Cohen-Macaulay$\,${\rm ;}
\item[{\rm (b)}] $\Omega_{A/k}$ is of linear type if and only if
${\rm edim}\, A_{\mathfrak p}\leq 2\dim A_{\mathfrak p}-1$ for
every non-minimal prime ${\mathfrak p}\in \spec (A)$$\,${\rm ;}
\item[{\rm (c)}] In case $A$ is normal, $\Omega_{A/k}$ is of linear type and normal
if and only if ${\rm edim}\, A_{\mathfrak p}\leq 2\dim A_{\mathfrak p}-2$ for
every non-regular prime ${\mathfrak p}\in \spec (A)$.
\end{enumerate}
\end{Proposition}
\demo Write $A=R/I$ as above.

(a) Since $\Omega_{A/k}$ satisfies $(F_0)$ it follows that $A$ is a complete
intersection locally in codimension one.
Furthermore, (\ref{conormalsequence}) implies $\depth I/I^2\geq \dim A-1$.
Finally notice that
\[(\wedge^2 I/I^2)^{**}\simeq \omega_A^* \simeq
\Hom_A(\omega_A\otimes _A \omega_A, \omega_A) \simeq \mbox{\rm Hom}_A(S_2(\omega_A),\omega_A).
\]
By \cite[1.3]{VaVi}, $S_2(\omega_A)$ is a maximal Cohen-Macaulay
$A$-module and therefore $(\wedge^2 I/I^2)^{**}$ is Cohen-Macaulay
as well. These facts imply that (\ref{Koszulcodim2}) is  exact by
the acyclicity lemma and that $\depth \za \geq \dim A+n-2$. On the
other hand, $\dim \za=\dim A+n-2$  since we are assuming that
$\Omega_{A/k}$ satisfies condition $(F_0)$ (see \cite[2.2]{sv}).
Therefore $\za$ is Cohen-Macaulay.

(b)  From part (a) we have that $\za$ is unmixed. In this case, the
module $\Omega_{A/k}$ is of linear type if and only if it satisfies
$(F_1)$  (see \cite[3.3 and the first remark on page 346]{sv}).
Alternatively, one could use the exactness of (\ref{Koszulcodim2}).

(c) If $\Omega_{A/k}$ satisfies $(F_2)$ then it is of linear type by part (b), i.e., $\za=\ra$.
By (a) $\za$ is Cohen--Macaulay. Thus, $\za=\ra$ satisfies condition (ii) of
Proposition~\ref{normalityandreflexiveness}, hence
the proposition implies that $\Omega_{A/k}$ is normal (and the
equality $\za=\ba$ holds).

Conversely, assume that $\Omega_{A/k}$ is of linear
type and normal.
For any non-regular prime ${\mathfrak p}\in\spec(A)$ write
$A_{\mathfrak p}=R/I$, where $(R,{\mathfrak n})$ is a regular local ring essentially
of finite type over $k$
 and $I\subset{\mathfrak n}^2$ is an $R$-ideal of height $g\leq 2$.
 By Theorem~\ref{normality_reflexiveness},
it suffices to show that $\mu(I+{\mathfrak n}^3/{\mathfrak n}^3)\leq 2g-1$.
This is clear if $g=1$ or else $g=2$ and $\mu(I)\leq 3$.
Thus we may assume that $g=2$ and $\mu(I)\geq 4$.
But in this situation the Hilbert--Burch theorem gives $I\subset {\mathfrak n}^3$.
\qed

\bigskip

We now treat the case where $A$ is a Gorenstein algebra of embedding codimension $3$. It has
a striking similarity to complete intersections. Write $A=R/I$, where $R$ is a regular
local ring
essentially of finite type over a perfect field $k$ and $I$ an ideal of height $3$.
As in (\ref{Koszulcodim2}), starting from a presentation of $\Omega_{A/k}$ we obtain
the ${\cal Z}$-complex
\begin{eqnarray}\label{Koszulcodim3}\nonumber
\kern-25pt 0 \lar (\wedge^3 I/I^2)^{**}\otimes_A B[-3] \lar
(\wedge^2 I/I^2)^{**}\otimes_A B[-2] &\lar &
( I/I^2)\otimes_A B[-1]
\lar B\\
 &\lar&  \mathbb{S}_{A/k}\lar 0\,. 
\end{eqnarray}

\begin{Proposition} \label{normcodim3}
Let $A$ be a reduced local Gorenstein ring essentially of finite type over a perfect
field $k$. Assume that $\ecodim A\leq 3$ and that $\Omega_{A/k}$ satisfies $(F_0)$.
\begin{enumerate}
\item[{\rm (a)}] The sequence {\rm (\ref{Koszulcodim3})} is  exact and
$\mathbb{S}_{A/k}$ is Gorenstein$\,${\rm ;}
\item[{\rm (b)}] $\Omega_{A/k}$ is of linear type if and only if
${\rm edim}\, A_{\mathfrak p}\leq 2\dim A_{\mathfrak p}-1$ for
every non-minimal prime ${\mathfrak p}\in \spec (A)$$\,${\rm ;}
\item[{\rm (c)}] In case $A$ is normal, $\Omega_{A/k}$ is of linear type and normal
if and only if ${\rm edim}\, A_{\mathfrak p}\leq 2\dim A_{\mathfrak p}-2$ for
every non-regular prime ${\mathfrak p}\in \spec (A)$.
\end{enumerate}
\end{Proposition}
\demo Since $\Omega_{A/k}$ satisfies $(F_0)$ it follows that $A$ is a complete
intersection locally in
codimension two.
We note that $(\wedge^3 I/I^2)^{**}$ is the determinant divisor of
 $I/I^2$, which is $A$ itself, and that the pairing
\[ \wedge^2 I/I^2 \times I/I^2 \rar \wedge^3 I/I^2 \rar A  \]
identifies $(\wedge^2 I/I^2)^{**}$ with $(I/I^2)^*$.
Furthermore $I/I^2$ is a Cohen-Macaulay $A$-module (see \cite[3.3(a)]{HeDefo}).
We conclude that the three left most modules
in (\ref{Koszulcodim3}) are maximal Cohen-Macaulay $B$-modules.
Now the argument proceeds as in the proof of Proposition~\ref{normcodim2}.
As for the Gorensteiness of $\mathbb{S}_{A/k}$ one uses the fact that (\ref{Koszulcodim3})
is a self-dual complex of $B$-modules.
\qed

\bigskip

The previous results motivate the following question:

\begin{Question} \rm Let $A$ be a local Cohen--Macaulay normal domain essentially of finite
type over a perfect field $k$. If $\Omega_{A/k}$ is
of linear type and normal does it follow that ${\rm edim}\,
A_{\mathfrak p}\leq 2\dim A_{\mathfrak p}-2$ for
every non-regular prime ${\mathfrak p}\in \spec (A)$ (or, equivalently,
${\mathbb R}_{A/k}={\mathbb B}_{A/k}$)?
\end{Question}
As to the converse, it was remarked earlier in the proof of
Proposition~\ref{normcodim2} that, in any codimension,
if $\za$ satisfies condition $(S_2)$ (e.g., if it is Cohen--Macaulay) then the inequalities
${\rm edim}\, A_{\mathfrak p}\leq 2\dim A_{\mathfrak p}-2$ for
every non-regular prime ${\mathfrak p}\in \spec (A)$ imply the normality of
$\Omega_{A/k}$ and also
the equality ${\mathbb R}_{A/k}={\mathbb B}_{A/k}$.

\bigskip

\bigskip

In the remainder of this section we will prove Theorem~\ref{normalisF2}, which is a
refined version of Proposition~\ref{normcodim2}(c).
For this we need several auxiliary results that may be of independent interest.
Recall that if $A$ is a Noetherian ring and $U\subset E$ are finite modules having a rank
then $U$ is said to be a {\it reduction\/} of $E$ if the induced inclusion
${\cal R}(U)\subset {\cal R}(E)$ is an integral ring extension (see, e.g., \cite{ram1}).

\begin{Proposition}\label{doubledualnonsense}  Let $A$ be a Noetherian ring
and  let $\phi: A^{g+1}\lar A^g$ be a homomorphism such that
$\mbox{\rm im}(\phi)$ has rank $g$ and $\Ext^1_A({\rm im}(\phi),
A)=0$. Then $I_g(\phi)\simeq \ker(\phi)^*$. In particular, if $A$
satisfies $(S_2)$ it follows that $I_g(\phi)$ is either the unit
ideal or an unmixed ideal of height one. If in addition $A$ is
Gorenstein locally in codimension one, then  ${\rm im}(\phi)$ is a
reduction of ${\rm im}(\phi)^{**}$.
\end{Proposition}
\demo
There is an exact sequence
\[  A^n \stackrel{\psi}{\lar} A^{g+1} \stackrel{\phi}{\lar} A^g,  \]
so that the entries of the first column of $\psi$ are the signed
maximal minors of $\phi$, hence generate $I_g(\phi)$, an ideal
of positive grade. Set $K= \ker (\phi)$ and $L={\rm im}(\phi)$. From 
the assumption, we obtain a short exact sequence
\[ 0 \rar L^{*} \lar {A^{g+1}}^* \lar K^* \rar 0\,,\]
which shows that $K^*\simeq {\rm im}(\psi^*)$.
Next we project ${A^n}^*$ onto the free module generated by the first
basis element, thus getting a commutative diagram
$$
\begin{array}{ccc}
{A^{g+1}}^* & \stackrel{\psi^*}{\lar} & {A^n}^*\\
\Vert && \downarrow\pi\\
{A^{g+1}}^* & \stackrel{\pi\psi^*}{\lar} & \kern-8pt A.
\end{array}
$$
Notice that ${\rm im}(\pi\psi^*)=I_g(\phi)$.
The projection $\pi$ induces a surjection ${\rm im}(\psi^*)\surjects
{\rm im}(\pi\psi^*)=I_g(\phi)$, which is necessarily an isomorphism since the
two modules are torsionfree of rank one.

Moreover, if $A$ satisfies $(S_2)$ then  $I_g(\phi)$ satisfies $(S_2)$,
hence it is the unit ideal or an unmixed ideal of height one.
 The last assertion is
a property of the integral closure of modules with
divisorial determinantal ideal (see \cite[the proof of 2.5]{icm}). \qed

\begin{Corollary}\label{acinonobstructed}
Let $R$ be a regular local ring essentially of finite type over a field.
Let $A=R/I$ be an almost complete intersection ring which is a
complete intersection locally in codimension one and satisfies $(S_2)$.
If $A$ is non-obstructed {\rm (}i.e., ${\rm T}^2(A/R,A)
={\rm Ext}^1_A(I/I^2,A)=0${\rm )}, then $I/I^2$ is a reduction of $(I/I^2)^{**}$.
\end{Corollary}
\demo Write $g=\height I$. Notice that $I/I^2$ has rank $g$ as an $A$-module
because $A$ is equidimensional.
Since $I$ is generated by $g+1$ elements, the first Koszul homology ${\rm H}_1(I)$
of these elements is the canonical module of $A$, hence satisfies $(S_2)$.
The exact sequence
$$0\rar {\rm H}_1(I)\rar A^{g+1}\rar I/I^2\rar 0$$
shows that $I/I^2$ is a torsionfree $A$-module.
Embedding $I/I^2$ into $A^g$ and applying Proposition~\ref{doubledualnonsense}
with ${\rm im}(\phi)=I/I^2$, we deduce the result.
 \qed

\begin{Proposition}\label{omeganotnormal}
Let $(A,\fm)$ be a normal local ring essentially of finite type over a perfect
field $k$ satisfying the following conditions$\,${\rm :}
\begin{enumerate}
\item[{\rm (a)}] $\dim A=\ecodim A=2\,${\rm ;}
\item[{\rm (b)}] $A$ is either a complete intersection or else an
 almost complete intersection defined by an ideal of order
$\geq 3$.
\end{enumerate}
Then $\Omega_{A/k}/\tau_{A}(\Omega_{A/k})$ is not integrally closed.
\end{Proposition}
\demo Write $A=R/I$,
where $(R,{\mathfrak n})$ is a regular local ring essentially of finite
type over $k$ and $I$ an ideal of height $2$. We may assume that $R$
has infinite residue field.
Set $\Omega=\Omega_{A/k}$.

In the complete intersection case, the
exact sequence
 $$0\rar I/I^2\rar A^n\rar \Omega\rar 0$$
 shows that $\Omega$ is a torsionfree $A$-module of projective
 dimension one, hence the assertion follows from Theorem~\ref{BurchGoto}.

Thus, we may assume that $A$ is an almost complete intersection and that $I\subset
{\mathfrak n}^3$.
There is an exact sequence
\begin{equation}\label{omegaone}
0\rar (I/I^2)^{**}\rar A^n\rar \Omega^1=\Omega/\tau_A(\Omega)\rar 0.
\end{equation}
By Theorem~\ref{Bergerlite}, $\Omega^1$ is not Cohen--Macaulay, hence not
reflexive. Thus $\fm\in \ass ((\Omega^1)^{**}/\Omega^1)$.
Suppose that $\Omega^1$ is integrally closed. In particular, it is $\fm$-full.
Thus, according to Proposition~\ref{mfullmod}, there exists an $x\in\fm$ such that
$A/\fm$ is a direct summand of $\Omega^1/x\Omega^1$.
Write  $(A',\fm')$ for the local ring $A/(x)$ and let $\{e_1,\ldots,e_n\}$
be the canonical basis of $A'\,^n$.
Tensoring (\ref{omegaone}) with $A'$ we obtain a presentation
\begin{equation}\label{omegaonebis}
0\rar M\rar A'\,^n \rar \Omega^1/x\Omega^1\rar 0,
\end{equation}
where we may assume that $M=\fm'\,e_1\oplus N$ with $N\subset \oplus_{i=2}^n A'e_i$.
Hence the image of $M$ under the projection $A'\,^n\surjects A'e_1\simeq A'$ is
$\fm'$.

On the other hand, as $A$ is a Cohen--Macaulay ring of embedding codimension $2$
it is non-obstructed (see \cite[3.2(a)]{HeDefo}).
Therefore by Corollary~\ref{acinonobstructed}, $I/I^2$ is a reduction of $(I/I^2)^{**}$.
Since $I\subset {\mathfrak n}^3$ it follows that $I/I^2\subset {\mathfrak n}^2A^n$.
Hence, projecting onto $A'e_1\simeq A'$ we see that $(\fm')^2$ is a reduction of $\fm'$.
This is impossible because $\dim A'>0$.
\qed

\begin{Theorem}\label{normalisF2}
Let $(A,\fm)$ be a normal local Cohen--Macaulay ring essentially of finite type over a perfect
field $k$ satisfying the following conditions$\,${\rm :}
\begin{enumerate}
\item[{\rm (a)}] $\ecodim A=2\,${\rm ;}
\item[{\rm (b)}] Locally in codimension $2$, the ring $A$ is either a complete
intersection or else an almost complete intersection defined by an ideal of order $\geq 3$.
\end{enumerate}
Then the following conditions are equivalent$\,${\rm :}
\begin{enumerate}
\item[{\rm (i)}] $\Omega_{A/k}$ is normal$\,${\rm ;}
\item[{\rm (ii)}] ${\rm edim}\, A_{\mathfrak p}\leq 2\dim A_{\mathfrak p}-2$ for
every non-regular prime ${\mathfrak p}\in \spec (A)$.
\end{enumerate}
Moreover, under any of the above equivalent conditions
$\Omega_{A/k}$ is of linear type and ${\mathbb R}_{A/k}$ is
Cohen--Macaulay.
\end{Theorem}
\demo  (i) $\Rightarrow$ (ii) Proposition~\ref{omeganotnormal} implies
that $\Omega_{A/k}$ satisfies $(F_1)$.
Hence by Proposition~\ref{normcodim2}(b), $\Omega_{A/k}$ is of linear type and
then according to Proposition~\ref{normcodim2}(c), it satisfies $(F_2)$.

(ii) $\Rightarrow$ (i) This follows from Proposition~\ref{normcodim2}(c).

The remaining assertions follow from the same proposition.
\qed

\subsection{Relation to Calabi--Yau varieties}

In this last part we explain to what extent the present results relate to
Calabi--Yau varieties.

Let $X\subset \pp^{n-1}_{\cc}$ be an arithmetically normal projective variety and
let $A$ stand for its
homogeneous coordinate ring. We say that $X$ is of {\it Calabi--Yau type\/}
if there
exists a homogenous isomorphism $\omega_A\simeq A$.
The notion of Calabi--Yau variety would also require that $X$ be
smooth and ${\rm H}^1(X,{\cal O}_X)=0$.
If $X$ is of Calabi--Yau type it often turns out that  $Y={\rm Proj}{(\ra)}
\subset \pp^{2n-1}_{\cc}$
has the same property.

We first look into the case of a complete intersection.

\begin{Proposition}
Let $X\subset \pp^{n-1}$ be a non-degenerate smooth projective variety that is a
complete intersection
of hypersurfaces of degrees $d_1\geq \cdots \geq d_g$, and let $A$ stand for its
homogeneous coordinate ring.
Assume that $X$ is of Calabi--Yau type and consider the subschemes $Y={\rm Proj}(\rac)
\subset Z={\rm Proj}(\zac) \subset \pp^{2n-1}_{\cc}$.
\begin{enumerate}
\item[{\rm (a)}] $Z$ is the complete intersection of $2g$ hypersurfaces of degrees
$d_1,\ldots,d_g, d_1,\ldots,d_g$ and $\omega_{\,{\zac}}\simeq
\zac$ as graded modules$\,${\rm ;}
\item[{\rm (b)}] If $d_1=2$ then $Z$ is neither reduced nor irreducible, and $Y$ is
not arithmetically Cohen--Macaulay$\,${\rm ;}
\item[{\rm (c)}] If $d_1\geq 3$ then $Y=Z$ is reduced and irreducible$\,${\rm ;}
\item[{\rm (d)}] The subscheme $Y$ is arithmetically normal if and
only if  $d_1\geq 4$ or else $d_2\geq 3$, in which
case $Y$ is of Calabi--Yau type. 
\end{enumerate}
\end{Proposition}
\demo Since $\omega_A\simeq A(-n+\sum_{i=1}^g d_i)$ and $X$ is
Calabi--Yau, we have $n=\sum_{i=1}^g d_i$. Since $X$ is
non-degenerate, $d_g\geq 2$. Therefore $n\geq 2g$, hence
$\Omega_{A/\cc}$ satisfies condition $(F_0)$. Since ${\rm
projdim}_A(\Omega_{A/\cc})=1$ this implies that $\zac$ is a complete
intersection (see \cite[Proposition 4]{Avramov}). Moreover, $\omega_{\,\,{\zac}}\simeq
\zac(-2n+2\sum_{i=1}^g d_i)=\zac$ as graded modules. This shows (a).

To prove (b), note that the hypothesis forces the equality $n=2g$ and
thus $\Omega_{A/\cc}$ does not satisfy condition $(F_1)$. It follows
from \cite[2.2]{BVV} that $Z$ is not irreducible, whereas
Theorem~\ref{redtorsionfree} implies that $Z$ is not reduced either.
Furthermore $Y$ is not arithmetically
Cohen--Macaulay according to Theorem~\ref{cmversusF1}.

As to (c), the assumption implies that $n\geq 2g+1$, hence
$\Omega_{A/\cc}$ satisfies condition $(F_1)$. Again since ${\rm
projdim}_A(\Omega_{A/\cc})=1$ we have $\zac=\rac$ (see
\cite[Proposition 4]{Avramov}, \cite[1.1]{Huneke}, \cite[3.4]{sv}).
Finally, notice that $d_1\geq 4$ or $d_2\geq 3$ if and only if
$n\geq 2g+2$. Thus, part (d) follows from
Corollary~\ref{CIandnormality}. \qed

\medskip

Next is the non-complete intersection case.
\begin{Proposition}
Let $X\subset \pp^{n-1}$ be a smooth arithmetically
Cohen--Macaulay projective variety that is not a complete
intersection. Assume that $\dim X\geq 2, \ecodim X\leq 3$ and $X$ is
of Calabi--Yau type. Consider the subschemes $Y={\rm Proj}(\rac)
\subset Z={\rm Proj}(\zac) \subset \pp^{2n-1}_{\cc}$, where $A$ stand
for the  homogeneous
coordinate ring of $X$ in the given embedding.
\begin{enumerate}
\item[{\rm (a)}] $Z$ is arithmetically Gorenstein and $\omega_{{\zac}}\simeq \zac$
as graded modules$\,${\rm ;}
\item[{\rm (b)}] If $\dim X=2$ then $Z$ is neither reduced nor irreducible$\,${\rm ;}
\item[{\rm (c)}] If $\dim X\geq 3$ then $Y=Z$ is reduced and irreducible$\,${\rm ;}
\item[{\rm (d)}] $Z$ is arithmetically normal if and only if $\dim X\geq 4$,
in which case $Y=Z$ is of Calabi--Yau type.
\end{enumerate}
\end{Proposition}
\demo Notice that $A= R/I$, where $R=\cc[X_1,\ldots,X_n]$ and
$I\subset (X_1,\ldots,X_n)^2$ is a homogeneous Gorenstein ideal of
height at most $3$, hence exactly of height $3$ because $X$ is not a
complete intersection. Since $X$ is smooth and $n=\hht\, I+\dim
A\geq 3+3=6$, the module $\Omega_{A/\cc}$ satisfies $(F_0)$. By
Proposition~\ref{normcodim3}(a), the ring $\zac$ is Gorenstein. In
fact, as a complex of graded modules over the standard graded
polynomial ring $B=\cc[X_1,\ldots,X_n,T_1,\ldots,T_n]$, the exact
sequence (\ref{Koszulcodim3}) now reads
$$
 0 \lar (\wedge^3 I/I^2)^{**}\otimes_A B \lar (\wedge^2
I/I^2)^{**}\otimes_A B \lar  ( I/I^2)\otimes_A B \lar B
 \lar  \mathbb{S}_{A/k}\lar 0\,,
$$
where $-^*={\rm Hom}_A(-,A)={\rm Hom}_A(-,\omega_A)$. From this we
get $\omega_{{\zac}}\simeq \zac$ as graded modules, proving part
(a).

As to (b), if $\dim X\geq 2$ then as above we see that
$\Omega_{A/\cc}$ does not satisfy $(F_1)$. Thus, $Z$ is not
irreducible according to \cite[2.2]{BVV}. On the other hand, by
considering the Hilbert function of $A$ modulo a linear system of
parameters one sees that $\dim _{\cc}[I]\,_2\leq 5=n-1$. Now an
application of Theorem~\ref{redtorsionfree} yields that $Z$ is not
reduced either, proving (b).

Finally, part (c) follows from Proposition~\ref{normcodim3}(b) and
part (d) from Proposition~\ref{normcodim3}(c). \qed

\bigskip

{\bf Acknowledgements:} The first author thanks the Department of
Mathematics 
at Purdue University
for hospitality during the preparation of part of this work. The other
two authors
 thank the Departmento de Matem\'atica at Universidade Federal
de Pernambuco (Brazil)
for hospitality during  various visits.

\end{document}